\documentclass[12pt, openany]{amsart}
\usepackage{import}
\usepackage[utf8]{inputenc}
\usepackage[margin=1.1in]{geometry}
\usepackage[italian,english]{babel}
\usepackage{csquotes}

\usepackage[maxbibnames=5, sorting=anyt, style=numeric, doi=false,isbn=false,url=false,eprint=false]{biblatex}
\usepackage{xcolor}
\usepackage{soulutf8}
\usepackage{appendix}
\usepackage[shortlabels]{enumitem} 
\usepackage{float} 
\usepackage{caption} 
\usepackage{tabularx}
\usepackage{makecell}
\usepackage{setspace}
\usepackage{amsmath,mathtools,amssymb,amsthm} %AMSMATH
\usepackage{stmaryrd,mathcomp,gensymb} %Symbols
\usepackage[makeroom]{cancel} %Canceling terms
\usepackage{tikz-cd}
\usetikzlibrary{arrows,fit,positioning,calc}
\usepackage{IEEEtrantools}

\usepackage{etoolbox} %Macro utili
\usepackage{cleveref}

\input{CustomTheorems.sty}
\input{CustomMath.sty}

\title{A theta operator for the group $\GSp_4$}
\author{Leonardo Fiore}
\date{August 2023}

\begin{document}

\maketitle

\begin{abstract}
    We construct a differential operator on sheaves of $p$-adic modular forms defined over the locus of $p$-rank $\ge 1$ of the Siegel threefold, by applying a revisited version of the approach that S.\ Howe recently introduced in \cite{howe} to construct the theta operator in the elliptic case.
\end{abstract}

\tableofcontents

\section{Introduction}

\subsection{Theta operators}
In 1977, Katz constructed a weight-raising differential operator on the space of $p$-adic modular forms, known as the $\vartheta$ operator. At the level of $q$-expansions, the operator $\vartheta$ can be described as $q\frac{d}{dq}$, meaning that $\vartheta(\sum_n a_n q^n)=\sum_n n a_n q^n$.

The construction that Katz presents in \cite{katz} goes as follows: if $X_1(N)\to \Spec(\ZZ[1/N])$ is the compactified modular curve of level $\Gamma_1(N)$ for some $N\ge 4$, and $\mathfrak{X}_1(N)\to \Spf(\Zp)$ denotes its $p$-adic completion at a prime $p\nmid N$, one has that, over the ordinary locus $\mathfrak{X}^{\ord}_1(N)\subset \mathfrak{X}_1(N)$, the Hodge filtration
$0 \to \omega\to H^1_{\mathrm{dR}} \to \omega^\vee\to 0$
splits canonically, via the so-called unit-root splitting, and this allows to derive a differential operator $\vartheta: \omega^k \to \omega^{k+2}$ (for any $k\in \ZZ$) from the Gauss-Manin connection on $H^1_{\mathrm{dR}}$. This operator has some interesting commutation properties with Hecke operators $T_\ell$ and $U_\ell$ at each prime $\ell$:
\begin{equation*}
    \begin{aligned}
        &T_\ell \vartheta = \ell T_\ell \vartheta& \text{ for }\ell\nmid N;\\
        &U_\ell \vartheta = \ell U_\ell \vartheta& \text{ for }\ell\mid N.\\
    \end{aligned}
\end{equation*}
As a consequence, applying the theta operator to a modular form corresponds to twisting by a cyclotomic character its associated Galois representation, as the following result formalizes (see, for example, the introduction of \cite{emfv} for a more detailed discussion).
\begin{theoremIntro}
    Let $f$ be a normalized cuspidal eigenform of weight $k\in \ZZ/(p-1)\ZZ$, level $\Gamma_1(N)$ and nebentypus $\epsilon: (\ZZ/N\ZZ)^\times \to E^\times$, with coefficients in some fine extension $E/\mathbb{F}_p$. Let $\rho_f: G_\QQ\to \GL_2(E)$ be the continuous semisimple Galois representation attached to it. If $g:=\vartheta(f)$, then $g$ is still a normalized cuspidal eigenform with coefficients in $E$ for the group $\Gamma_1(N)$, whose weight is $k+2$, and $\rho_{g} \cong \rho_f\otimes \omega$, being $\omega: G_\QQ\to \GL_1(E)$ the cyclotomic character.
\end{theoremIntro}

Recent progress has been done in reproducing Katz's construction to more general Shimura varieties, and in the consequent study of the resulting differential operators, both in the mod $p$ and in the $p$-adic setting. In 2012, E.\ Eishen published her work \cite{e1}, which discusses how to apply Katz's construction to Shimura varieties attached to unitary groups of the form $U(n,n)$. A generalization of the result appeared in the 2018 paper \cite{emfv} by Eishen, Fintzen, Mantovan and Varma, which addresses more general Shimura varieties of type A and C having non-empty ordinary locus.

When the ordinary locus is empty, one can still adopt a definition similar to Katz's one on the $\mu$-ordinary locus, after suitably replacing the unit-root splitting, which is no longer available, with a subtler canonical splitting of $H^1_{\mathrm{dR}}$, associated to the slope filtration of the universal $\mu$-ordinary $p$-divisible group. This approach was explored by De Shalit and Goren in \cite{dg}, and by Eischen and Mantovan in \cite{em}.

A quite different perspective on $\vartheta$ operators was proposed in the paper \cite{howe} by S.\ Howe, appeared in 2020. The author restricts its attention to the elliptic case, and starts by considering the big Igusa tower on the $p$-adic modular curve $\mathfrak{X}_1(N)\to\Spf(\Zp)$, i.e.\ the moduli space $\mathfrak{IG}\to \mathfrak{X}_1(N)$ classifying the trivialisations $A[p^\infty]\cong \mupi \oplus \QpZp$ of the universal $p$-divisible group $A[p^\infty]$. On $\mathfrak{IG}$, the group of self-$p$-quasi-isogenies of $\mupi \oplus \QpZp$ acts, which yields a residual infinitesimal action by the $p$-divisible group $\mupi$ on Katz's space of modular forms. By taking the derivative of this action, one recovers Katz's $\vartheta$ operator.

Howe's approach has several advantages. Its geometric nature, together with its independence from $q$-expantion computations, makes is well-suited for being abstracted and applied to more general Shimura varieties of PEL type. Moreover, the fact that the theta operator is obtained as the differential of a $\mupi$-action ensures that it can be iterated $p$-adically, without any need of proving this via explicit congruences between powers of $\vartheta$.

\subsection{The Siegel threefold}

In the present work, we will focus our attention exclusively on the Siegel threefold, which is to say the Shimura variety attached to the reductive group $\GSp_4$, that can be seen as the moduli of principally polarized abelian surfaces (with some prime-to-$p$ level structure, that we will always consider fixed, and for which we presume that the neatness assumption is satisfied). Let $X$ a toroidal compactification of the Siegel threefold, $A\to X$ be the semiabilian surface it carries, and $\omega$ the sheaf of invariant differentials of $A$. We denote by $\mathfrak{X}\to \Spf(\Zp)$ the $p$-adic completion of $X$.

Over the ordinary locus $\mathfrak{X}^{\ord}\subset \mathfrak{X}$, theta operators can easily be constructed by using a mild generalisation of Howe's approach, with the sole difference that the infinitesimal group acting turns now out to be $\mupi^3$, and not just $\mupi$ as in the elliptic case, so that differentiation is now possible in three distinct directions. Equivalently, one can more classically construct theta operators by following \cite{emfv}'s perspective, which is to say adopting Katz's original approach and exploiting the unit-root splitting of $H^1_{\mathrm{dR}}$.

We will not further examine this case, as the focus of the present paper is studying theta operators over loci of $\mathfrak{X}$ larger than the ordinary one, for which the generalization of Howe's approach is still possible, but more involved. In his 2020 work \cite{pilloniGSP4}, devising the foundations of higher Hida and Coleman theories in the case of $\GSp_4$, Pilloni is concerned with defining and studying spaces of $p$-adic modular forms defined over the $p$-rank $\ge 1$ locus, $\mathfrak{X}^{\ge 1}$. This is larger that the ordinary locus (which can be described as the $p$-rank $2$ locus), it is still an open formal subscheme of $\mathfrak{X}$, and it is not affine: it is, instead, the union of two open formal affine subschemes. The non-affineness of $\mathfrak{X}^{\ge 1}$ makes it interesting, for example, for $p$-adically deforming those Hecke eigenclasses that lie in the $H^1$ groups of modular sheaves -- which would vanish if restricted to the ordinary locus. The construction of $p$-adic modular forms on $\mathfrak{X}^{\ge 1}$ proposed in \cite{pilloniGSP4} goes as follows. The first step is constructing an affine pro-étale cover $\mathfrak{X}^{\ge 1}_{\mathrm{Kli}(p^\infty)}\to \mathfrak{X}^{\ge 1}$ classifying multiplicative height-1 subgroups $H\inj A$, being $A$ the semi-abelian surface carried by $\mathfrak{X}$. Over this deep Klingen cover, he defines a (small) Igusa variety $\mathfrak{Ig}$ classyfing the trivializations $H\cong \mupi$, and he exploits to define, for each $p$-adic character $k$ of $\Zp^\times$, a sheaf of $p$-adic modular forms $\mathfrak{F}^k$ as the sheaf of functions on $\mathfrak{Ig}$ on which $\Zp^\times$ acts via the character $k$.

Given $k_1$ a $p$-adic character of $\Zp^\times$ and $k_2\in \ZZ$, we will adopt the notation $\mathfrak{F}^{(k_1, k_2)}$ to denote the sheaf $\mathfrak{F}^{k_1-k_2}\otimes \det(\omega)^{k_2}$.

\subsection{Our main results}
Retaining the notations introduced in the previous subsection, the main result of this paper is the construction of a differential operator $\vartheta: \mathfrak{F}^{(k_1, k_2)}\to \mathfrak{F}^{(k_1+2, k_2)}$. More precisely, we have the following result, which is \Cref{thm main result}.
\begin{theoremIntro}[A]
        For any pair $(k_1,k_2)$ with $k_1: \Zp^\times\to R$ a character (for $R$ some $p$-adically complete $\Zp$-algebra) and $k_2\in \ZZ$, each multiplicative continuous function $f: \Zp\to R$ induces an $R$-linear operator $\vartheta^f: \mathfrak{F}^{(k_1,k_2)}\to \mathfrak{F}^{(k_1+2f,k_2)}$ over $\mathfrak{X}_{\mathrm{Kli}(p^\infty)}^{\ge 1}\cotimes R$ minus the boundary. In particular, for $n\in \NN$, the function $\Zp \to \Zp$, $x\mapsto x^n$ induces a $n$-th order $R$-linear differential operator $\vartheta^n: \mathfrak{F}^{(k_1,k_2)}\to \mathfrak{F}^{(k_1+2n,k_2)}$.
\end{theoremIntro}

Following a variant of Howe's approach, we deduce this operator by first constructing an action by the infinitesimal group $\mupi$ on an appropriate (small) Igusa covering of $\mathfrak{X}_{\mathrm{Kli}}(p^\infty)$, that we denote $\mathfrak{T}^\gr$ in the paper and is a $(\Zp^\times \times \GL_1)$-torsor over $\mathfrak{X}_{\mathrm{Kli}(p^\infty)}$, and by then taking the differential of this action. We build the $\mupi$-action by attaching to each element $\zeta\in \mupi$ of the acting group a deformation of the universal $p$-divisible group $A[p^\infty]$, and then applying Serre-Tate lifting theory -- without the need to consider big Igusa towers and quasi-isogeny groups acting on them, as in Howe's original manner.

This newly constructed $\vartheta$ operator has suitable commutation properties with Hecke operators (see \Cref{prop commutation with Hecke prime to p}), which allow us to connect its action with the twist by the cyclotomic character on the Galois representation side, as expressed in \cref{thm Galois}, whose statement we report here for convenience.
\begin{theoremIntro}[B]
    Let $f\in H^i(X,\omega^{(k_1,k_2)})\otimes \Qpbar$ be an eigenform for the Hecke algebra $\mathcal{H}:=C_c^0(\GSp_4(\mathbb{A}_f)//K)$, being $i=0$ or $i=1$; let $\rho_f$ be its attached Galois representation (see \cite[Theorem 5.3.1]{pilloniGSP4}). Let $g:=\vartheta(f)\in H^i(\mathfrak{Y}^{\ge 1}_{\mathrm{Kli}(p^\infty)},\mathfrak{F}^{(k_1+2,k_2)})\otimes \Qpbar$, where $\mathfrak{Y}^{\ge 1}_{\mathrm{Kli}(p^\infty)}$ denotes $\mathfrak{X}^{\ge 1}_{\mathrm{Kli}(p^\infty)}$ minus the boundary. Assume $g\neq 0$. Then, $g$ is still an eigenform for the prime-to-$p$ Hecke algebra $\mathcal{H}^{p}$; moreover, there exists a Galois representation $\rho_g$ attached to its Hecke eigensystem, which coincides with the cyclotomic twist of $f$ (in other words, $\rho_g\cong \rho_f\otimes \omega$, being $\omega$ the cyclotomic character).
\end{theoremIntro}
\subsection{Outline}

The paper is organised as follows. \Cref{sec preliminaries} recalls some useful algebra results; \Cref{sec setup} discusses the structure of the Siegel threefold and sets up the moduli spaces of interest for the construction. The actual construction is performed in \Cref{sec constructing}. \Cref{sec weight hecke} is devoted to studying the effect $\vartheta$ has on weights and Hecke eigenvalues. Finally, \Cref{sec applications} links the newly constructed $\vartheta$ operator to cyclotomic twists of Galois representations.
\subsection{Acknowledgments}
The author is sincerely grateful to his PhD advisor Fabrizio Andreatta for proposing him the problem, and, most importantly, for the continuous and patient guidance, suggestions and insights that were essential to this research work.

\section{Preliminaries and conventions}
\label{sec preliminaries}
\subsection{Actions}
All group actions will be taken to be left actions. Given a base scheme $S$, an action of an $S$-group scheme $G$ on an $S$-scheme $\pi: X\to S$ is by definition a morphism of $S$-schemes $\rho: G\times_S X \to X$ satisfying the action axioms. Every section $g\in G(S)$ induces an automorphism $\psi_g: X\to X$. The action of $G$ on $X$ also induces an action of $G$ on the sheaf $\pi_*\OO_X$ of the functions on $X$ (and, more in general, on the cohomology $R\pi_*\OO_X$). A section $s\in \pi_*(\OO_X)(S)$ can be viewed as a scheme homomorphism  $s: X\to \mathbb{A}^1_S$, and the element $g\in G(S)$ acts by taking the $s$ to $s \circ \psi_g^{-1}$.

\subsection{Vector bundles}
Given a base scheme $S$, a rank-$d$ vector bundle $V$ on $S$ can alternatively be defined as
\begin{enumerate}[(a)]
    \item a finite locally-free sheaf of $\OO_S$-modules of rank $d$, 
    \item a scheme over $S$ that is locally isomorphic to an affine space $\mathbb{A}^d_S$ (where ``locally'' can equivalently mean locally for the Zariski, étale, fppf, or fpqc topology), or
    \item a $GL_d$-torsor (for the Zariski, the étale, the fppf, or the fpqc topology, equivalently).
\end{enumerate}

The definitions (a) and (b) correspond in the following way. If $V$ is a vector bundle in the sense (a), the corresponding affine scheme over $S$ is obtained as $\Spec_S(\Sym_{\OO_S}(V^\vee))$. Instead, if $E$ is a vector bundle in the sense (b), then the corresponding $\OO_S$-sheaf is the sheaf of sections of the structure map $\pi: E\to S$.  %Let us show that the two operations are inverse to each other. If $V$ is a vector bundle in the sense (a), then a section of $\Spec_S(\Sym_{\OO_S}(V^\vee))$ over $S$ is just an $\OO_S$-algebra morphism $\Sym_{\OO_S}(V^\vee)\to \OO_S$, hence an $\OO_S$-linear morphism $V^\vee\to \OO_S$, and hence an section of $V$ over $S$. Conversely, if $\pi: E\to S$ is a vector bundle in the sense (b), let $\mathcal{F}$ of sections of $\pi: E\to S$. The morphism $\Spec_S(\Sym_{\OO_S}(\mathcal{F}^\vee))\to E$ is given by taking 

We now turn to the relation between (a) and (c). If $V$ is a finite locally-free $\OO_S$-module of rank $d$, the corresponding left $\GL_d$-torsor is $\mathcal{T} = \Isom(V, \OO_S^d)$, on which $g\in \GL_d$ acts via $\psi\mapsto g\circ \psi$. Conversely, if $\mathcal{T}$ is a $\GL_d$-torsor, the sheaf of $\GL_d$-invariant morphisms of $S$-schemes $\mathcal{T}\to \mathbb{A}^d_S$ gives the corresponding finite, locally free $\OO_S$-module of rank $d$.

\subsection{Torsors and representations}

Let $G$ be a group scheme over a base scheme $S$, and let $\pi: \mathcal{T}\to S$ be a left $G$-torsor over $S$. For every left representation $V$ of $G$ over $S$ (i.e., $V$ is a finite locally free $\OO_S$-module, endowed with an action $\rho: G\to \GL(V)$), one can form the $S$-scheme $\mathcal{T}^\rho:=\mathcal{T}\times^G V=\{(x,v): x\in \mathcal{T}, v\in V\}/\{(x,v)\sim(gx,gv): g\in G\}$; in other words, $\mathcal{T}^\rho$ is the quotient of the $S$-scheme $\mathcal{T}\times_S V$ under the diagonal action of $G$. 
\begin{proposition}
    The $S$-scheme $\mathcal{T}^\rho$ is a vector bundle; as a finite, locally-free $\OO_S$-module, it corresponds to the sheaf of $G$-equivariant morphisms of $S$-schemes $\mathcal{T}\to V$.
    \begin{proof}
        Let $\psi: G\isom \mathcal{T}$ be a trivialization of the torsor $\mathcal{T}$ over some cover $S'\to S$. Then, $f_\psi: V\to \mathcal{T}^\rho, v\mapsto (\psi(1_G), v)$ and $g_\psi: \mathcal{T}^\rho\to V$, $(x,v)\mapsto (\psi^{-1}(x) v)$ are clearly inverse to each other, and hence $\mathcal{T}^\rho$ gets identified with $V$ over $S'$. In particular, $\mathcal{T}^\rho$ is a vector bundle. Given a $G$-equivariant morphism $h: \mathcal{T}\to V$ over $S'$, then the element $(\psi(1_G),h(\psi(1_G)))\in \mathcal{T}^\rho(S')$ is clearly independent of the trivialization $\psi$ chosen, and hence $h$ canonically corresponds to a section of $\mathcal{T}^\rho$ over $S'$. It is easy to see that, conversely, every section of $\mathcal{T}^\rho$ determines a unique equivariant morphism $\mathcal{T}\to V$.
    \end{proof}
\end{proposition}

The proof of the proposition emphasises that $\mathcal{T}^\rho$ can be thought of as a $\mathcal{T}$-twisted version of $V$, in the sense that any given trivialisation of $\mathcal{T}$ induces an isomorphism $V\cong\mathcal{T}^\rho$. This also emerges very clearly from the two examples given below.
\begin{example}
    \label{example rep line bundle}
    Let $\rho: G\to \GM{S}$ be a character of the group $G$. In this case, $\mathcal{T}^\rho$ is a line bundle, and it coincides with $(\pi_*\OO_{\mathcal{T}})[-\rho]$, i.e.\ the invertible $\OO_S$-subsheaf of $\pi_*\OO_{\mathcal{T}}$ consisting of those functions on $\mathcal{T}$ on which $G$ acts via its character $-\rho$.
\end{example}
\begin{example}
    \label{example GL2 rep vec bundle}
    Suppose $G=\GL_2$, let $W$ be a rank-2 vector bundle on $S$, and let $\mathcal{T}$ be the corresponding left $G$-torsor $\Isom(W,\OO_S^2)$. Suppose we are given a dominant weight $(k_1,k_2)\in \ZZ^2$ for $G$ (the dominance condition is $k_1\ge k_2$), and let $V_{k_1,k_2}$ be the irreducible representation of highest weight $(k_1, k_2)$ of $\GL_2$, i.e.\  $V_{k_1,k_2}=\Sym^{k_1-k_2}(\OO_S^2)\otimes \det^{k_2}(\OO_S^2)$ with the obvious canonical left action of $\GL_2=\Aut(\OO_S^2)$. Then, we have that $\mathcal{T}^\rho \cong \Sym^{k_1-k_2}(W)\otimes \det^{k_2}(W)$.
\end{example}

\subsection{$p$-divisible groups} We denote by $\Nilp_{\Zp}$ the site of all $\Zp$-algebras on which $p$ is nilpotent. Given $R\in \Nilp(R)$, the stack of all $p$-divisible groups (also known as Barsotti-Tate groups) on $R$ is denoted $\BT(R)$. Given an abelian scheme $A$ over $R$, its attached $p$-divisible group is denoted $A[p^\infty]$. Any prime-to-$p$ polarization $\lambda: A\to A^\vee$ (i.e., a positive symmetric isogeny of degree prime to $p$) induces an anti-symmetric isomorphism $\lambda: A[p^\infty]\isom A[p^\infty]^\vee$. This justifies the following definition.
\begin{definition}
    \label{def polarized bt group}
    A \emph{polarized} $p$-divisible group on $R$ is a $p$-divisible group $G\in \BT(R)$ endowed with an isomorphism $\lambda: G\isom G^\vee$ satisfying $\lambda^\vee = -\lambda$.
\end{definition}
\begin{remark}
    This notion of polarized $p$-divisible group agrees with that of ``$(-1)$-polarized $p$-divisible group'' in \cite[\S3.1.1]{moonen}.
\end{remark}

Fix now $I\le R$ a nilpotent ideal. Suppose we are given a polarized $p$-divisible group $G_0\in \BT(R/I)$, and a $p$-divisible group $G\in \BT(R)$ deforming $G_0$. The polarization $\lambda_0: G_0\to G_0^\vee$ may or may not lift to a polarization $\lambda: G\to G^\vee$; in case it lifts, the lifting $\lambda$ is necessarily unique by Drinfeld's rigidity lemma (see \cite[Lemma 1.1.3]{KatzSerreTate}), and we say that $G$ is a polarized deformation of $G_0$, or a deformation of $G_0$ as a polarized $p$-divisible group. Serre-Tate lifting theory ensures that deforming the (polarized) $p$-divisible group of a principally polarized abelian scheme is equivalent to deforming the whole abelian scheme.
\begin{theorem}
    \label{thm serre tate}
    Fix a ring $R\in \Nilp_{\Zp}$, a nilpotent ideal $I\le R$, and a principally polarized abelian scheme $A_0$ over $R_0:=R/I$, with $p$-divisible group $G_0:=A_0[p^\infty]$. For every deformation $(G, \psi_G: G\otimes R/I\isom G_0)$ of the polarized $p$-divisible group $G_0$ to a polarized $p$-divisible group $G$ on $R$, there exists one and only one deformation $(A, \psi_A: A\otimes R/I\isom A_0)$ of $A_0$ to a principally polarized abelian scheme on $R$, such that the $p$-divisible group of $A$ coincides with $G$, meaning that there exists a (necessarily unique) isomorphism $\varphi: A[p^\infty]\isom G$ such that, modulo $I$, the following diagram commutes:
    \begin{equation*}
    \begin{tikzcd}
        A[p^\infty]\otimes R/I\ar["\psi_A"]{d}\ar["\varphi"]{r}& G\otimes R/I \ar["\psi_G"]{d}\\
        A_0[p^\infty]\ar[equal]{r}&G_0
    \end{tikzcd}
    \end{equation*}
    \begin{proof}
        It follows from Serre-Tate lifting theorem (see \cite[Theorem 1.2.1]{KatzSerreTate}) that the abelian scheme $A$ deforming $A_0$ whose $p$-divisible group coincides with $G$ exists unique. From the same theorem, it follows that there exists a unique isogeny $\lambda: A\to A^\vee$ that lifts the principal polarization on $A_0$. Now, to verify that $\lambda$ is a principal polarization one has to prove that it is symmetric, positive and an isomorphism, but all three conditions can be verified fiberwise, and hence they are satisfied since they are already known to hold modulo $I$.
    \end{proof}
\end{theorem}

%\subsection{Dominant weights and subgroups}
%Suppose $G$ is a reductive group over a field $F$; let $P$ a parabolic subgroup, with Levi quotient $L=P/U$. If $V$ is a representation of $G$, 

\subsection{The Iwasawa algebra of a profinite group}
In this subsection, we recall some introductory notions about measures and distributions on profinite groups; our reference is \cite{mazur}.

Let $R$ be a $p$-adically complete $\Zp$-algebra. Given a profinite group $G$, one can define its group algebra $R[[G]]:=\varprojlim_H R[G/H]$, where $H$ ranges among all open normal subgroups of $H$, and $G/H$ consequently ranges among all finite discrete quotients of $H$. This is also known as the Iwasawa algebra of $G$; it is a commutative formal Hopf algebra (where sums and products are computed pointwise), and its elements may be thought of as $R$-valued \emph{measures} on the group $G$, where a $R$-valued measure on $G$ is a rule $\mu$ attaching to each closed and open subset $E$ of $G$ a scalar $\mu(E)\in R$, in such a way that $\mu(\emptyset)=0$ and $\mu(E_1\sqcup \ldots \sqcup E_n)=\mu(E_1)+\ldots+\mu(E_n)$.

The Iwasawa algebra $R[[G]]$ is dual to the co-commutative formal Hopf algebra of continuous functions $\Cont(G,R)$, the duality pairing being given by integration: $(\mu, f)\mapsto \int f d\mu$, for each $\mu\in R[[G]]$, $f\in \Cont(G,R)$.

\subsection{Actions by $\mupi$}
\label{subsec action by continuous functions}
In this subsection, we will present equivalent ways of describing an action by the infinitesimal group $\mupi$ on a formal scheme.
\begin{proposition}
    Let $R$ be a $p$-adically complete $\Zp$-algebra. The formal group $\mupi\to \Spf(R)$ is the formal spectrum of the Iwasawa algebra of the profinite group $\Zp$. In other words, there exists a canonical isomorphism of formal group schemes $\mupi\cong \Spf(R[[\Zp]])$.
    \begin{proof}
        Given $A\in \Nilp_R$, an $A$-point of $\Spf(R[[\Zp]])$ is a continuous algebra homomorphism $R[[\Zp]]\to A$, which necessarily factors through some $R[\ZZ/p^n\ZZ]$ giving rise to a morphism $R[\ZZ/p^n\ZZ]\to A$, i.e.\ a group homomorphism $\ZZ/p^n\ZZ\to A^\times$, which is the same as a $p^n$-root of unity $\zeta$ in $A$, i.e.\ an $A$-point of $\mupi$.
    \end{proof}
\end{proposition}

\begin{proposition}
    \label{prop equivzlent action mupi}
    Suppose $\mathfrak{X}$ is a formal group scheme over some $p$-adically complete base ring $R$. Then, the following data are equivalent:
    \begin{enumerate}[(a)]
        \item a left action by $\mupi$ on $\mathfrak{X}$;
        \item a right co-action by the formal Hopf algebra $R[[\Zp]]$ on the sheaf of algebras $\OO_{\mathfrak{X}}$;
        \item a left action by the formal Hopf algebra $\Cont(\Zp,R)$ on the sheaf of algebras $\OO_{\mathfrak{X}}$.
    \end{enumerate}
    \begin{proof}
        We will only show how a left $\mupi$-actions yields the algebra actions (b) and (c). Since the group $\mupi$ is infinitesimal, it acts on $\mathfrak{X}$ via universal homeomorphisms; in particular, the action preserves all open subsets of the formal scheme $\mathfrak{X}$; if $f$ is a function defined on some open subscheme $U$ of $\mathfrak{X}$, and $\zeta\in \mupi(R)$, then one can define $\zeta\cdot f$ to be $f: U\to \mathbb{A}^1$ precomposed with the automorphism $U\to U$ induced by $\zeta^{-1}$. This induces a left action by $\mupi$ on $\OO_{\mathfrak{X}}$, which is to say a right co-module structure on $\OO_{\mathfrak{X}}$ for the Hopf algebra $R[[\Zp]]$, given by a morphism $\OO_{\mathfrak{X}}\to \OO_{\mathfrak{X}}\cotimes R[[\Zp]]$. But this is the same as a left module structure on $\OO_{\mathfrak{X}}$ by the dual formal Hopf algebra $\Cont(\Zp, R)$.
    \end{proof}
\end{proposition}

Given a section $s$ of $\OO_{\mathfrak{X}}$ over an open affine subset $U\subseteq \mathfrak{X}$, an element $\zeta\in \mupi(R)$ (i.e., an element $\zeta\in R^\times$ such that $\zeta^{p^n}\to 1$ as $n\to \infty$): we want to construct $\zeta \cdot s$, and this can be done in two equivalent ways.
\begin{enumerate}
    \item we see $s$ as a function $U\to \mathbb{A}^1$, and precompose it with the automorphism $\psi_{\zeta^{-1}}: U\to U$ given by the action of $\zeta\in \mupi(R)$: in other words, $\zeta\cdot s := s \circ \psi_{\zeta^{-1}}$;
    \item we consider the continuous function $\Zp \to R$, $x\mapsto \zeta^x$, which is the \emph{group-like} element of $\Cont(\Zp,R)$ corresponding to $\zeta$, and we let it act on $s$ adopting the point of view of \Cref{prop equivzlent action mupi}(c); in other words, $\zeta\cdot s := [x\mapsto \zeta^x] \cdot s$.
\end{enumerate}
One particular section of the infinitesimal group $\mupi$ is given by $\zeta=1+\varepsilon\in \mupi(R[\varepsilon])$, where $\epsilon^2=0$. For each function $f\in \OO_{\mathfrak{X}}(U)$, we have that $(1+\varepsilon) f = f + \varepsilon f'$, for some section $f'\in\OO_\mathfrak{X}(U)$. As $(1+\varepsilon)$ is the canonical base element for the 1-dimensional Lie algebra of $\mupi$, the operation $f\mapsto f'$ can be thought of as the \emph{derivative} of the action of $\mupi$ on $\mathfrak{X}$. The continuous character $\Zp \to R[\varepsilon]$ corresponding to $\zeta=1+\varepsilon$ is $x\mapsto (1+\varepsilon)^x = 1+x\varepsilon$, which shows that $f'$ can alternatively be defined as the image of $f$ through the action of the continuous function $\Zp\to R$, $x\mapsto x$. Note that $[x\mapsto x]$ is a primitive element of the formal Hopf algebra $\Cont(\Zp, R)$, which entails that the Leibnitz rule holds: $(f_1 f_2)'= f_1'f_2+f_1f_2'$.
\begin{definition}
    \label{def induced theta}
    Let $\mathfrak{X}$ be a formal scheme over a $p$-adically complete $\Zp$-algebra $R$, carrying an $R$-linear $\mupi$-action on $\mathfrak{X}$. Then, the \emph{differential of the action} is the $R$-linear differential operator $\OO_{\mathfrak{X}}\to \OO_{\mathfrak{X}}, f\mapsto f'$ introduced above. 
\end{definition}

\section{The $\GSp_4$ setting}
\label{sec setup}

In this section, we present the construction of $p$-adic modular forms on $\mathfrak{X}^{\ge 1}$, following, for almost all the material presented, Pilloni's paper \cite{pilloniGSP4}.

\subsection{The Siegel threefold}
Let $G=\GSp_4$, and fix a prime $p$. Let $X\to \Spec(\ZZ_p)$ denote a toroidal compactification of the Siegel threefold of neat level $K\le G$, where $K=K_pK^p$ with $K_p=\GSp_4(\ZZ_p)$. The subgroup $K_\ell$ equals the maximal compact subgroup $\GSp_4(\ZZ_\ell)$ for all but finitely many primes; let $N$ be the product of those finitely many primes (so that no level structure is present away from $N$). The threefold $X$ carries a canonical semi-abelian surface $A\to X$. We denote by $D$  the boundary divisor, and we write $Y:=X\setminus D$ for the non-compactified threefold, i.e.\ the open subscheme of $X$ over which $A$ is an abelian surface (actually, $A_{|Y}\to Y$ is the \emph{universal} principally polarized abelian surface of level $K$). 

The sheaf of invariant differentials of $A$, which we denote $\omega$, is a rank-2 vector bundle on $X$, i.e.\ a torsor for the group $M:=\GL_2$, which should be thought of as the standard Levi subgroup of $G=\GSp_4$, modulo its center. For each choice of $(k_1,k_2)\in \ZZ^2$, $k_1\ge k_2$, we define the automorphic vector bundle $\omega^{(k_1,k_2)}:=\Sym^{k_1-k_2}(\omega)\otimes \det^{k_2}(\omega)$ corresponding to the irreducible representation of highest weight $(k_1,k_2)$ of the group $M$ (compare with \Cref{example GL2 rep vec bundle}).

\subsection{Hecke algebra}
\label{subsec hecke alg}
The cohomology spaces $H^i(X,\omega^{(k_1, k_2)})\otimes \Qp$ of the automorphic vector bundles are finite-dimensional vector spaces carrying an action of the Hecke algebra $\mathcal{H}=C_c^0(\GSp_4(\mathbb{A}_f)//K)$, where $C_c^0$ denotes the continuous, compactly supported functions. At each prime $\ell\nmid N$, the Hecke algebra $\mathcal{H}_\ell$ is $\mathbb{Z}[T_{\ell,0}, T_{\ell,1}, T_{\ell,2}, T_{\ell,0}^{-1}]$, being $T_{\ell,0}, T_{\ell,1}$ and $T_{\ell,2}$ the characteristic functions of the double cosets $\diag(\ell,\ell,\ell,\ell)$, $\diag(\ell^2,\ell,\ell,1)$, and $\diag(\ell,\ell,1,1)$. Each of these Hecke operators $T_{\ell,i}$ admits an interpretation as a correspondence from $X\otimes \Qp$ to itself; more precisely,
\begin{enumerate}[(a)]
    \item the correspondence inducing $T_{\ell,0}$ classifies the isogenies $A\to A'$ whose kernel is the subgroup $A[\ell]$;
    \item the correspondence inducing $T_{\ell,1}$ classifies the isogenies $A\to A'$ whose kernel is a totally isotropic subgroup of $A[\ell^2]$, whose intersection with $A[\ell]$ has degree $\ell^3$.
    \item the correspondence inducing $T_{\ell,2}$ classifies the isogenies $A\to A'$ whose kernel is a totally isotropic subgroup of $A[\ell]$.
\end{enumerate}
To summarize, the operator $T_{\ell,i}$ is associated with isogenies whose kernel is a totally isotropic subgroup of $A[\ell^{M(i)}]$, where $M(0)=M(1)=2,$ and $M(2)=1$. For $\ell\neq p$, these correspondences are immediately seen to be also defined integrally. For $\ell=p$, the integral definition is subtler, and the integrality properties of the operators are discussed in \cite[\S 7]{pilloniGSP4}.

For later use, we also introduce here the Hecke polynomial $Q_\ell(X)\in \mathcal{H}_\ell[X]$, which is given by
    \begin{equation*}
        Q_\ell(X)=1
        -T_{\ell,2} X
        +\ell(T_{\ell,1}+(\ell^2+1)T_{\ell,0}) X^2
        - \ell^3 T_{\ell,2} T_{\ell,0} X^3
        + \ell^6 T^2_{\ell,0} X^4.
    \end{equation*}
    \begin{remark}
    \label{rmk Ql homogeneous}
    The Hecke polynomial is homogeneous of degree 0, if one assigns degree $-M(i)$ to the symbol $T_{\ell,i}$ and degree $1$ to the symbol $X$.
    \end{remark}

\subsection{Klingen level}
Let $\mathfrak{X}\to \Spf(\ZZ_p)$ be the formal completion of $X$ at $p$, and let $\mathfrak{Y}$ be that of $Y$. We will work on the open formal loci $\mathfrak{X}^{\ge 1}$ and $\mathfrak{Y}^{\ge 1}$ where the $p$-rank of $A$ is at least one.

The locus $\mathfrak{X}^{\ge 1}$  is stable under the action of Hecke operators; contary to the ordinary locus, it is not affine; instead, it is the union of two distinct, formal affine subschemes. This can be verified via the theory of partial Hasse invariants: the ordinary locus (i.e., the locus $\mathfrak{X}^{=2}$ where the $p$-rank is 2) is the complement, in $\mathfrak{X}$, of the zero-locus $\mathfrak{X}^{\le 1}$ of the Hasse invariant $\mathrm{Ha}\in H^0(\mathfrak{X}\otimes \Fp, \det^{p-1}(\omega))$. Moreover, a second Hasse invariant $\mathrm{Ha}'\in H^0(\mathfrak{X}^{\le 1}\otimes \Fp, \det^{p^2-1}(\omega))$ can be defined so that it vanishes exactly over the $p$-rank $0$ locus: see \cite[\S 6.3]{pilloniGSP4} for its definition. Now, since $\det(\omega)$ is an ample line bundle on $\mathfrak{X}\otimes \Fp$, a power $(\mathrm{Ha}')^k$ necessarily lifts (non-canonically) to a section $\widetilde{(\mathrm{Ha}')^k}\in H^0(\mathfrak{X}\otimes \Fp, \det^{k(p^2-1)}(\omega))$. Recalling that the complement of an ample divisor on a proper scheme is always affine, one can now conclude that $\mathfrak{X}^{\ge 1}=\lbrace \mathrm{Ha}\neq 0\rbrace \cup \lbrace \widetilde{(\mathrm{Ha}')^k}\neq 0\rbrace$ is the union of two affines.

For each $n\ge 1$, we can define a cover $\mathfrak{X}_{\mathrm{Kli}(p^n)}^{\ge 1}\to \mathfrak{X}^{\ge 1}$ classifying the subgroups $H_n\le A[p^n]$ that are étale-locally isomorphic to $\mu_{p^n}$. The map $\mathfrak{X}_{\mathrm{Kli}(p^n)}^{\ge 1}\to \mathfrak{X}^{\ge 1}$ is affine and étale (see \cite[Lemma 9.1.1.1]{pilloniGSP4}), but it is not finite: over the ordinary locus $\mathfrak{X}^{=2}\subset \mathfrak{X}^{\ge 1}$, it has rank $p+1$, while it is an isomorphism when restricted to the closed subscheme $\mathfrak{X}^{=1}\subset \mathfrak{X}^{\ge 1}$ where the $p$-rank is 1. By taking the inverse limit, one obtains a pro-étale cover $\mathfrak{X}_{\mathrm{Kli}(p^\infty)}^{\ge 1}\to \mathfrak{X}^{\ge 1}$, which classifies all subgroups $H\le A[p^\infty]$ that are pro-étale-locally isomorphic to $\mu_{p^\infty}$.

On $\mathfrak{Y}_{\mathrm{Kli}(p^\infty)}^{\ge 1}$, the inclusion of the universal multiplicative height-1 subgroup $H\inj A[p^\infty]$ together with the dual projection $A[p^\infty]\onto H^\vee$ induce a filtration $0=\Fil_0(A[p^\infty])\subset \Fil_1(A[p^\infty])\subset \Fil_2(A[p^\infty])\subset \Fil_3(A[p^\infty])=A[p^\infty]$; the graded piece $\gr_2(A[p^\infty])=H^\vee$ is étale of height 1, the graded piece $\gr_0(A[p^\infty])=H$ is multiplicative of height 1, while the middle graded piece $\gr_1(A[p^\infty])$ is a $p$-divisible group of height 2 and dimension 1, which can alternatively be ordinary or connected-connected, depending on whether the $p$-rank of $A$ is 1 or 2 (here, by $\gr_i(A[p^\infty])$ we mean $\Fil_{i+1}(A[p^\infty])/\Fil_{i}(A[p^\infty])$). At the level of differentials, we get an induced filtration  $0=\Fil^3(\omega)=\Fil^2(\omega)\subset \Fil^1(\omega)\subset \Fil^0(\omega)=\omega$, where $\gr^i(\omega):=\Fil^{i}(\omega)/\Fil^{i+1}(\omega)$ is the sheaf of invariant differentials of the $p$-divisible group $\gr_i(A[p^\infty])$. In particular, $\gr^2(\omega)=0$, $\gr^1(\omega)=\ker(\omega\onto\omega_H)$, $\gr^0(\omega)=\omega_H$. While the filtration $\Fil_i(A[p^\infty])$ only exists away from the boundary $D$, the filtration on the rank-2 vector bundle $\omega$ extends to the whole toroidal compactification $\mathfrak{X}^{\ge 1}_{\mathrm{Kli}}(p^\infty)$.

%Over $\mathfrak{X}_{\mathrm{Kli}(p^\infty)}^{\ge 1}$, the inclusion $H\inj A[p^\infty]$ induces a projection $\omega\onto \omega_H$, where $\omega_H$ is the line bundle of invariant differentials of $H$; hence, $\omega$ carries a filtration $0=\Fil_1(\omega)\subset \Fil_2(\omega)\subset \Fil_3(\omega)=\omega$, whose graded pieces are the line bundles $\gr^1(\omega)=\ker(\omega\onto \omega_H)$ and $\gr^0(\omega)=\omega_H$, where $\gr_i(\omega):=\Fil_{i+1}(\omega)/\Fil_{i}(\omega)$.

If we let $M^{\gr}\le M$ be the subgroup $\GL_1\times \GL_1\le \GL_2$, over $\mathfrak{X}_{\mathrm{Kli}(p^\infty)}^{\ge 1}$ we have the left $M$-torsor $\mathcal{T}$ of the trivializations $\psi: \omega\isom \OO^2$, and the left $M^{\gr}$-torsor $\mathcal{T}^{\gr}$ of the trivializations  $\psi_0\oplus \psi_1: \gr^0(\omega)\oplus \gr^1(\omega)\isom \OO\oplus \OO$.
For each choice of $(k_1,k_2)\in \ZZ^2$ with $k_1\ge k_2$, we have
\begin{enumerate}
    \item a dominant weight $(k_1, k_2)$ for $M$, corresponding to an irreducible representation $V_{k_1, k_2}$ of $M$ that, when twisted by the $M$-torsor $\mathcal{T}$, gives rise to the vector bundle $\omega^{(k_1, k_2)}=\Sym^{k_1-k_2}(\omega)\otimes \det^{k_2}(\omega)$ on  $\mathfrak{X}_{\mathrm{Kli}(p^\infty)}^{\ge 1}$; and
    \item a character $(k_1, k_2)$  for $M^{\gr}$, i.e.\ a 1-dimensional representation $V^{\gr}_{k_1,k_2}$ of $M^{\mathrm{\gr}}$, which, once twisted by $\mathcal{T}^\gr$, gives rise to the line bundle $\gr^0(\omega)^{k_1-k_2}\otimes \det^{k_2}(\omega)\cong \gr^0(\omega)^{k_1}\otimes \gr^1(\omega)^{k_2}$, which we will denote by $\mathfrak{F}^{(k_1,k_2)}$.
\end{enumerate}
The inclusion $M^{\gr}\le M$ gives a surjective equivariant morphism $V_{(k_1,k_2)}\onto V_{(k_1,k_2)}^{\gr}$, which in turn provides the canonical projection $\omega^{(k_1,k_2)}\onto \mathfrak{F}^{(k_1,k_2)}$. This is an isomorphism only for scalar weights (i.e., when $k_1=k_2$).

\subsection{$p$-adic interpolation}
\label{subsec p-adic interpolation}
Let now $\mathfrak{T}^{\gr}$ be the moduli space that classifies isomorphisms $\iota: \mupi\isom \gr_0 A[p^\infty]$ and $\psi_1: \gr^1(\omega)\isom \OO$ over $\mathfrak{X}_{\mathrm{Kli}(p^\infty)}^{\ge 1}$: this is a $\mathfrak{M}^{\gr}$-torsor, where $\mathfrak{M}^{\gr}\le M^\gr$ is the subgroup $\mathbb{Z}_p^\times \times \GL_1\le \GL_1\times \GL_1$. An element $(a_0, a_1)\in \mathfrak{M}^\gr$ acts by taking $(\iota,\psi_1)\mapsto (\iota\circ a_0, a_1\circ \psi_1)$. If we fix the standard trivialization $\OO\cong \omega_{\mupi}, 1\mapsto dt/t$, then the differential of the map $\iota$ clearly induces an isomorphism $d\iota: \gr^0(\omega)\isom \OO$, and the assignement $(\iota, \psi_1)\mapsto (d\iota, \psi_1)$ consequently defines an equivariant morphism from the $\mathfrak{M}^\gr$-torsor $\mathfrak{T}^\gr$ to the $M^{\gr}$-torsor $\mathcal{T}^\gr$, which identifies $\mathcal{T}^\gr$ as the pushforward of $\mathfrak{T}^\gr$ along the inclusion $\mathfrak{M}^\gr\le M$. In other words, the torsor $\mathfrak{T}^\gr$ we have just constructed is a $\mathfrak{M}^\gr$-reduction of the $M^\gr$-torsor $\mathcal{T}^\gr$.

This reduction of the structure group allows for a (partial) $p$-adic variation of the weight $(k_1,k_2)$. More precisely, given a continuous character $k_1: \Zp^\times\to R$ for $R$ a $p$-adically complete $\Zp$-algebra, and an integer $k_2\in \ZZ$, we can form the continuous character $(k_1, k_2): \mathfrak{M}^{\gr}\to R$, which gives rise to a line bundle  $(\mathfrak{T}^{\gr})^{(k_1, k_2)}$ over $\mathfrak{X}_{\mathrm{Kli}(p^\infty)}^{\ge 1}\cotimes_{\Zp} R$ (see \Cref{example rep line bundle}), which we will still denote by $\mathfrak{F}^{(k_1,k_2)}$ (when $k_1\in \mathbb{Z}$, this definition is compatible with the one given in the previous subsection). We will set $\mathfrak{F}^{k}:=\mathfrak{F}^{(k,0)}$, so that  $\mathfrak{F}^{(k_1,k_2)}=\mathfrak{F}^{k_1-k_2} \otimes \det(\omega)^{k_2} = \mathfrak{F}^{k_1}\otimes \gr^1(\omega)^{k_2}$. Finally, if we denote by $\kappa: \Zp^\times\to \Lambda:=\Zp[[\Zp^\times]]$ the universal character of $\Zp^\times$, we have that any continuous $p$-adic character $k$ of $\Zp^\times$ is given by a continuous algebra homomorphism $k: \Lambda\to R$ with values in some $p$-adically complete $\Zp$-algebra $R$, and $\mathfrak{F}^k=\mathfrak{F}^\kappa\cotimes_{\Lambda,k} R$. The sheaf $\mathfrak{F}^{\kappa}$ coincides with the one that Pilloni constructs in \cite[\S 9.3]{pilloniGSP4} and denotes the same way.

\section{Constructing the theta operator}
\label{sec constructing}

This section aims at constructing the $\vartheta$ operator on $\mathfrak{T}^\gr$, which is the $(\Zp^\times \times \GL_1)$-torsor on $\mathfrak{X}^{\ge 1}_{\mathrm{Kli}(p^\infty)}$ defined in the previous one. Following Howe's approach, we first construct a $\mupi$-action on $\mathfrak{T}^\gr$, and then take its differential to get the $\vartheta$ operator. We are actually only able to achieve these results away from the boundary $D$.

The actual construction is carried out in \Cref{subsec defining action}, after some necessary preliminaries on extensions of $p$-divisible groups presented in \Cref{subsec kummer}.

\subsection{Kummer extentions}
\label{subsec kummer}
Let $R\in \Nilp_{\Zp}$ be a $\Zp$-algebra on which $p$ is nilpotent, and let $\Ext^1(\QpZp, \mupi)(R)$ be the groupoid of the extensions of $\QpZp$ by $\mupi$ over $R$. An object of $\Ext^1(\QpZp, \mupi)(R)$ consists of a triple $(E, \iota, \pi)$, where $E$ is a $p$-divisible group, while $\iota: \mupi\inj E$ and $\pi: E\onto\QpZp$ are two morphisms such that $0\to \mupi\to E\to \QpZp\to 0$ is an exact sequence.  

The groupoid $\Ext^1(\QpZp, \mupi)(R)$ also carries a $\Zp$-linear structure:
\begin{enumerate}
    \item given two extensions $E_1$ and $E_2$, one can form their Baer sum $E_1+E_2$;
    \item given an extension $E\in \Ext^1(\QpZp, \mupi)(R)$ and a scalar $a\in \underline{\Zp}(R)$, it is possible to define the extension $aE$ as the pushout of $E$ along the morphism $\mupi\to \mupi$, $\zeta\mapsto \zeta^a$, or equivalently as the pullback of $E$ along the morphism $\QpZp\to \QpZp$, $m\mapsto am$: it is immediate to verify that the two results are canonically isomorphic.
\end{enumerate}
In particular, the group $\underline{\Zp^\times}(R)$ acts on $\Ext^1(\QpZp, \mupi)(R)$, and, for all $a\in \underline{\Zp^\times}(R)$, we have that $a\cdot (E, \iota, \pi)=(E,\iota, a^{-1}\pi) \cong (E, \iota a^{-1}, \pi)$.

The trivial extension $\mupi\oplus \QpZp$ is an object of $\Ext^1(\QpZp,\mupi)$; its automorphism group is $\Hom_R(\QpZp,\mupi)=T_p\mupi(R)$.

\begin{proposition}
    Every extension $E\in \Ext^1(\QpZp, \mupi)(R)$ becomes trivial over some pro-finite faithfully flat $R$-algebra $R'$. In other words, the map $\mathrm{triv}: \star\to \Ext^1(\QpZp,\mupi)$ is fpqc-surjective, and it induces an equivalence between $\Ext^1(\QpZp, \mupi)$ and the stack of torsors under the group scheme $T_p\mupi$.
    \begin{proof}
        Consider the element $p^{-n}\in \QpZp$, which generates the cyclic subgroup $\QpZp[p^n]\cong \mathbb{Z}/p^n \mathbb{Z}$. Since $E[p^n]\onto \QpZp[p^n]$ is a surjective map of finite flat schemes, $p^{-n}$ lifts to a section of $E[p^n]$ after a finite flat extention of the base ring $R$. By repeating this consideration for all $n$, one can form, after a pro-finite flat extension of $R$, a compatible system of lifts of all elements $p^{-n}\in \QpZp$, i.e.\ a splitting of $E\onto \QpZp$
    \end{proof}
\end{proposition}

Given an element $\zeta\in \mupi(R)$, one can interpret $\zeta$ as an $R$-linear morphism $\Zp\to \mupi, m\mapsto \zeta^m$ over $R$, and form an extension $E_\zeta\in\Ext^1(\QpZp,\mupi)(R)$ by pushing out the exact sequence $0\to \Zp\to \Qp\to \QpZp\to 0$ along $\zeta: \Zp\to \mupi$. This clearly defines a $\Zp$-linear map $\mupi(R)\to \Ext^1(\QpZp,\mupi)$. If we denote by $\Ext^1(\QpZp, \mupi)(R)^\tors$ the sub-groupoid consisting of all those extensions $E\in \Ext^1(\QpZp, \mupi)(R)$ such that $p^n E$ is trivial for some $n$, we have the following result.
\begin{proposition}
    \label{prop torsion ext equivalent conditions}
    For an object $(E, \iota, \pi)\in \Ext^1(\QpZp,\mupi)(R)$, the following are equivalent.
    \begin{enumerate}
        \item[(a)] $(E, \iota, \pi)$ is in the image of the map $\mupi(R)\to \Ext^1(\QpZp,\mupi)(R), \zeta\mapsto E_\zeta$ described above.
        \item[(b)] $(E, \iota, \pi)\in \Ext^1(\QpZp,\mupi)(R)^\tors$.
        %\item[(c)] $(E, \iota, \pi)\in \Ext^1(\QpZp,\mupi)(R)$ becomes trivial over $R/I$ for some nilpotent ideal $I\le R$.
        \item[(c)] There exists $q: E\to\mupi$ such that $q \iota = p^n$ for some $n$.
        \item[(d)] There exists $j: \QpZp\to E$ such that $\pi j = p^n$ for some $n$.

    \end{enumerate}
    \begin{proof}
        Since $\mupi(R)$ is a $p$-torsion group, it is completely clear that (a) implies (b). Let us now prove that (b) implies (c). The extension $p^n E$ is the pushout of $E$ along the morphism $p^n: \mupi\to \mupi$:
        \begin{equation*}
        \begin{tikzcd}
            0 \ar{r}& \mupi \ar["\iota"]{r}\ar["p^n"]{d} & E \ar["\pi"]{r}\ar{d} & \QpZp\ar{r}\ar[equal]{d} & 0\\
            0 \ar{r}& \mupi \ar["\iota'"]{r} & p^n E \ar["\pi'"]{r} & \QpZp\ar{r}& 0
        \end{tikzcd}
        \end{equation*}
        The existence of a splitting $q': p^n E\to \mupi$ of $\iota'$ clearly entails the existence of a morphism $q: E\to \mupi$ such that $q\iota=p^n$, which proves (c).
        That (c) and (d) are equivalent is completely clear. Finally, let us show that (d) implies (a). Let $j: \QpZp\to E$ be a morphism such that $\pi j = p^n$. Then, $p^{-n}j: \Qp\to E$ induces the following morphisms of extensions
        \begin{equation*}
        \begin{tikzcd}
            0 \ar{r}& \Zp \ar{r}\ar{d} & \Qp \ar{r}\ar["p^{-n}j"]{d} & \QpZp\ar{r}\ar[equal]{d} & 0\\
            0 \ar{r}& \mupi \ar{r} & E \ar{r} & \QpZp\ar{r}& 0
        \end{tikzcd}
        \end{equation*}
        which shows that $E=E_\zeta$ for some $\zeta\in \mupi(R)$.
    \end{proof}
\end{proposition}
\begin{remark}
    If $\zeta\in \mupi(R)$, then let $I:=(\zeta-1)R$, so that $E_\zeta$ is identified with the trivial extention $E_1=\mupi \oplus \QpZp$ modulo $I$. Since $\zeta$ is a $p$-th power root of unity in $R$, and $p$ is nilpotent in $R$, it is immediate to realize that the ideal $I$ is nilpotent. This proves that, for $\zeta\in \mupi(R)$, $E_\zeta\in \Ext^1(\QpZp,\mupi)(R)$ comes with a canonical trivialization modulo the nilpotent ideal $I:=(\zeta-1)R$; in other words, it is a \emph{deformation} of the trivial extension $\mupi\oplus \QpZp\in \Ext^1(\QpZp,\mupi)(R/I)$ over $R$.
\end{remark}
\begin{lemma}
    \label{lemma kummer}
    For each $\zeta\in \mupi(R)$, there exists a canonical polarization $E_\zeta\cong E_{\zeta}^{\vee}$, which is the unique polarization that, modulo the nilpotent ideal $(\zeta-1)$, coincides with the canonical polarization of $\mupi\oplus \QpZp$.
    \begin{proof}
        Omitted.
    \end{proof}
\end{lemma}

The following lemma about morphisms between Kummer extensions will also be useful later.
\begin{lemma}
    \label{lemma morphism Kummer}
    Given $R\in \Nilp_{\Zp}$, take $\zeta_1, \zeta_2\in \mupi(R)$, and take $a, b\in \Zp^\times$. Modulo the nilpotent ideal $I:=(\zeta_1-1, \zeta_2-1)$, both Kummer extensions $E_{\zeta_1}$ and $E_{\zeta_2}$ are canonically identified with $\mupi\oplus \QpZp$. If $(\zeta_1)^a = (\zeta_2)^b$, the morphism  $a\oplus b: \mupi\oplus \QpZp\to \mupi\oplus \QpZp$ over $R/I$ lifts to a (necessarily unique) morphism $g_{a,b}: E_{\zeta_1}\to E_{\zeta_2}$ over $R$. Such morphism restricts to $a$ on $\mupi$, and coincides with $b$ on the quotient $\QpZp$.
    \begin{proof}
        Since $E_{\zeta_i}$ is the pushout of $0\to \Zp\to \Qp\to \QpZp$ along $\zeta_i: \Zp\to \mupi$, and $\zeta_2=\zeta_1^{a/b}$, by functoriality there exists a unique morphism $h: E_{\zeta_1}\to E_{\zeta_2}$ that restricts to $ab^{-1}$ on the subgroup $\mupi$ and to the identity on the quotient $\QpZp$. If we compose $h$ with the $b$-homothety morphism $b: E_{\zeta_2}\to E_{\zeta_2}$, we get the desired lift $g_{a,b}$.
    \end{proof}
\end{lemma}

\subsection{Defining the action}
\label{subsec defining action}
Let $R\in \Nilp_{\ZZ_p}$, and let us consider a point $x: \Spec(R)\to \mathfrak{T}^{\gr}$, away from the boundary, where $\mathfrak{T}^{\gr}\to \mathfrak{X}_{\mathrm{Kli}(p^\infty)}^{\ge 1}$ is the $\mathfrak{M}^{\gr}$-torsor torsor defined in \S\ref{subsec p-adic interpolation}. Let $\Ab{x}$ be the corresponding abelian surface, and $\omega_x$ be its sheaf of invariant differential. We recall that the $p$-divisible group $\Ab{x}$ comes with a 3-step filtration
\begin{equation*}
    0 = 
    \Fil_0 (A(x)[p^\infty]) \subset 
    \Fil_1 (A(x)[p^\infty]) \subset 
    \Fil_2 (A(x)[p^\infty]) \subset 
    \Fil_3 (A(x)[p^\infty]) 
    = A(x)[p^\infty];
\end{equation*}
moreover, for the first and the last graded piece we have trivializations $\iota_x: \mupi \isom \gr_0 (\Ab{x}[p^\infty])$ and $\pi_x: \gr_2 (\Ab{x}[p^\infty])\isom \QpZp$ that are dual to each other with respect to the principal polarization carried by $\Ab{x}$; the middle piece $\gr_1 (\Ab{x}[p^\infty])$, instead, is not trivialized (actually, its isomorphism type varies), but its sheaf of invariant differentials $\gr^1(\omega_x)$ is.

Let now take $\zeta\in \mupi(R)$. We have shown in \S\ref{subsec kummer} how $\zeta$ gives rise to an extension $E_\zeta\in \Ext^1(\QpZp,\mupi)(R)$ that becomes canonically trivial modulo the nilpotent ideal $I=(\zeta-1)R$. We are now going to define a deformed version $E_\zeta+\Ab{x}[p^\infty]$ of the filtered $p$-divisible group $\Ab{x}[p^\infty]$; the definition can be given in two equivalent ways:
\begin{enumerate}
    \item we take the pushout $E_1$ of $E_\zeta$ along $\mupi\isom \gr_0(\Ab{x}[p^\infty])\inj \Fil_2(\Ab{x}[p^\infty])$; we have that $E_\zeta+\Ab{x}[p^\infty]$ can now be defined as the Baer sum of $E_1$ and $\Ab{x}[p^\infty]$ inside $\Ext^1(\Ab{x}[p^\infty]/\Fil_2(\Ab{x}[p^\infty]),\Fil_2(\Ab{x}[p^\infty]))$;
    \item we take the pullback $E_2$ of $E_\zeta$ along $\Ab{x}[p^\infty]/\Fil_1(\Ab{x}[p^\infty])\onto \gr_2(\Ab{x}[p^\infty])\isom \QpZp$; now, $E_\zeta+\Ab{x}[p^\infty]$ can be defined as the Baer sum of $E_2$ and $\Ab{x}[p^\infty]$ inside $\Ext^1(\Ab{x}[p^\infty]/\Fil_1(\Ab{x}[p^\infty]),\Fil_1(\Ab{x}[p^\infty]))$.
\end{enumerate}
\begin{lemma}
    The two results, obtained via procedures (1) and (2) above, coincide. Moreover, the resulting $p$-divisible group $E_\zeta+\Ab{x}[p^\infty]$ also comes with a 3-step filtration $\Fil_i(E_\zeta+\Ab{x}[p^\infty])$. 
    \begin{proof}
        See \cite[\S9.3]{ExtPan}.
    \end{proof}
\end{lemma}
\begin{remark}
    \label{rmk baer sum panache}
    By construction, the quotients $\Fil_i/\Fil_{i+d}$ for $E_\zeta+\Ab{x}[p^\infty]$ are canonically identified with those of  $\Ab{x}[p^\infty]$, for $0\le d\le 2$ (the difference between the two $p$-divisible groups only shows up for $d=3$). 
    In particular, the trivializations of $\gr^0(\Ab{x}[p^\infty])$, of $\gr^2(\Ab{x}[p^\infty])$ and of the sheaf of invariant differentials of $\gr^1(\Ab{x}[p^\infty])$ induce corresponding trivializations for $\Ab{x}[p^\infty]+E_\zeta$.
\end{remark} 

\begin{lemma}
    \label{BaerSumLemma}
    Let $G$ denote the filtered $p$-divisible group $A(x)[p^\infty]$. The $G\mapsto E_\zeta+G$ operation just introduced satisfies the following properties:
    \begin{enumerate}[(a)]
        \item $E_{\zeta_1}+(E_{\zeta_2}+G) = (E_{\zeta_1}+E_{\zeta_2})+G = E_{\zeta_1 \zeta_2}+G$.
        \item $(E_{\zeta}+G)^\vee \cong (E_{\zeta})^\vee+G^\vee$.
        \item If $\zeta=1$, $E_{\zeta}+G$ canonically coincides with $G$.
    \end{enumerate}
    \begin{proof}
        Omitted.
    \end{proof}
\end{lemma}

\begin{proposition}
    \label{prop def p-div group}
        The filtered $p$-divisible group $E_\zeta+\Ab{x}[p^\infty]$ is a deformation over $R$ of the filtered $p$-divisible group $\Ab{x}[p^\infty]\in \BT(R/I)$, being $I$ the nilpotent ideal $(\zeta-1)R$. Moreover, the polarization on $\Ab{x}[p^\infty]\in \BT(R/I)$ lifts (in a necessarily unique way) to a polarization on $E_\zeta+\Ab{x}[p^\infty]\in \BT(R)$.    
    \begin{proof}
        By \Cref{BaerSumLemma}(c), the filtered $p$-divisible group $E_{\zeta}+\Ab{x}[p^\infty]$ is a deformation of $\Ab{x}[p^\infty]\in \BT(R/I)$ over $R$. Moreover,  \Cref{BaerSumLemma}(b) ensures that $(E_{\zeta}+\Ab{x}[p^\infty])^\vee$ canonically coincides with $E_{\zeta}^\vee+\Ab{x}[p^\infty]^\vee$. Now, if we pack together the isomorphism $A(x)[p^\infty]\isom (A(x)[p^\infty])^\vee$ induced by the polarization on $A(x)$, together with the canonical identification $E_{\zeta}\isom (E_{\zeta})^\vee$ given by \Cref{lemma kummer}, one obtains a polarization $E_\zeta+\Ab{x}[p^\infty]\isom (E_\zeta+\Ab{x}[p^\infty])^\vee$ lifting the one carried by $\Ab{x}[p^\infty]\in \BT(R/I)$.
    \end{proof}
\end{proposition}

  We are now finally to construct the $\mupi$-action.
\begin{theorem}
    \label{prop existence theta action}
    There is a canonical action of the formal group $\mupi$ on the moduli space $\mathfrak{T}^{\gr}\cotimes R$ minus the boundary. Given $x\in \mathfrak{T}^{\gr}(R)$, corresponding to an abelian surface $\Ab{x}/R$, the $p$-divisible group of the abelian surface $\Ab{\zeta\cdot x}$ is the deformation $E_\zeta + \Ab{x}[p^\infty]$ constructed at the beginning of this subsection.
    \begin{proof}
        If has been noted in \Cref{prop def p-div group} that $E_\zeta+\Ab{x}[p^\infty]$ is a polarized deformation of $\Ab{x}[p^\infty]$. By Serre-Tate lifting theory (see \Cref{thm serre tate}), it follows from \Cref{prop def p-div group} that $A(x)/(R/I)$ deforms uniquely to a principally polarized abelian surface $\zeta \cdot A$ over $R$ whose $p$-divisbile group $(\zeta \cdot A)[p^\infty]$ is $E_\zeta+\Ab{x}[p^\infty]$ of $\Ab{x}[p^\infty]\in \BT(R/I)$. The polarized abelian surface $\zeta\cdot A$, together with the additional structure on its $p$-divisible group presented in \Cref{rmk baer sum panache}, is classified by a point $\zeta\cdot x: \Spec(R)\to\mathfrak{T}^{\gr}$. This assignment satisfies the action axioms thanks to \Cref{BaerSumLemma}(a) and \Cref{BaerSumLemma}(c).
    \end{proof}
\end{theorem}

By taking the differential of the action, we now define our $\vartheta$ operator.
\begin{definition}
    \label{prop theta perspective}
    The $\vartheta$ operator on $\mathfrak{T}^\gr\setminus D$ is the differential operator obtained by taking the differential of the the $\mupi$-action on $\mathfrak{T}^\gr\setminus D$ constructed in \cref{prop existence theta action}, in the sense of \Cref{def induced theta}.
\end{definition}

\section{Interaction with weight and Hecke operators}
\label{sec weight hecke}
On the scheme $\mathfrak{T}^{\gr}\setminus D$, along with the action of $\mupi$ we have constructed in \cref{prop existence theta action}, there is an action by $\mathfrak{M}^{\gr}=\mathbb{Z}_p^\times\times \GL_1$; while the latter respects the structure morphism $\pi: \mathfrak{T}^{\gr}\to \mathfrak{X}_{\mathrm{Kli}(p^\infty)}^{\ge 1}$, we remark that the first one does not (actually, it does not even respect the structure morphism $\mathfrak{T}^{\gr}\to \mathfrak{X}^{\ge 1}$). Moreover, one can also define Hecke operators acting on the cohomology on $\mathfrak{T}^\gr$, that we will make explicit below. The aim of this subsection is understanding how the action of $\mupi$ interacts with the action of $\mathfrak{M}^{\gr}$ (hence, with the \emph{weight} of modular forms) and with Hecke operators (hence, with the \emph{Hecke eigenvalues} of modular forms). 

\subsection{The effect on weights}
To begin with, we introduce some notation. Let $R\in \Nilp_{\Zp}$, let $x_1, x_2\in \mathfrak{T}^{\gr}(R)$ be two points away from the boundary, and suppose we are given an isogeny $f: \Ab{x_1}\to \Ab{x_2}$ between the abelian surfaces classified by those points. Assume that $f$ respects the filtrations on $\Ab{x_1}[p^\infty]$ and $\Ab{x_2}[p^\infty]$, so that $f$ restricts to a morphism $\gr_i f$ between the graded pieces $\gr_i(\Ab{x_1}[p^\infty])\to \gr_i(\Ab{x_2}[p^\infty])$ for all $i=0,1,2$. Both the domain and the codomain of $\gr_2 f$ are identified with $\QpZp$, hence $\gr_2 f$ is just the multiplication by some scalar $m_2(f)\in \underline{\Zp}(R)$. Analogously, the domain and the codomain of $\gr_0 f$ are idenfied with $\mupi$, hence $\gr_0 f$ is a power map $\mupi\to \mupi$, $\zeta\mapsto \zeta^{m_0(f)}$ for some $m_0(f)\in \underline{\Zp}(R)$. Finally, since the sheaf of differentials of the $p$-divisible groups $\gr_1(\Ab{x_i}[p^\infty])$ are trivizlied, the map $d(\gr_1(f)): \gr_1(\omega_{x_2})\to \gr_1 (\omega_{x_1})$ is just $R\to R$, $a\mapsto m_1(f) a$ for some scalar $m_1(f)\in R$.

\begin{lemma}
    \label{rmk m_i for id}
    The following hold.
    \begin{enumerate}[(a)]
        \item     Let $x\in \mathfrak{M}^{\gr}(R)$ and $(a_0, a_1)\in \mathfrak{M}^\gr(R)$. Then, the identity morphism $f: \Ab{x}\to \Ab{(a_0, a_1)\cdot x}$ satisfies $m_2(f)=a_0, m_1(f)=a_1^{-1}$ and $m_0(f)=a_0^{-1}$.
        \item Let $x_1, x_2, x_3\in  \mathfrak{M}^{\gr}(R)$, and let $f: A(x_1)\to A(x_2)$ and $g: A(x_2)\to A(x_3)$ be isogenies respecting the filtrations on $p$-divisible groups. Then, $m_i(gf) = m_i(g) m_i(f)$ for all $i=1,2,3$.
        \item If $\lambda: A(x)\isom A(x)^\vee$ is the principal polarization of which $A(x)$ is endowed, then $m_i(\lambda)=1$ for all $i=1,2,3$.
        \item If $N: A(x)\to A(x)$ is the multiplication-by-$N$ isogeny, then $m_i(N)=N$ for all $i=1,2,3$.
        \item Let $f: A(x_1)\to A(x_2)$ be an isogeny whose kernel is a totally isotropic subgroup of $A(x_1)[N]$, for some $N\ge 1$. Suppose that $f$ respect the polarizations $\lambda_i$ carried by $A(x_i)$, meaning that $f^\vee \lambda_2 f = N\lambda_1$. Then, $m_2(f)m_0(f)=N$
    \end{enumerate}
    \begin{proof}
        We omit the verification of (a), (b), (c), and (d). To prove (e), we notice that we have $f^\vee \lambda_2 f = N\lambda_1$ by assumption. If we apply $m_0$ to the left-hand side one obtains $m_0(f^\vee) m_0(\lambda_2) m_0(f) = m_2(f) \cdot 1 \cdot m_0(f)$. By applying it on the right-hand side, one gets $N m_0(\lambda_2) = N$, whence  the result follows.
    \end{proof}
\end{lemma}

\begin{lemma}
    \label{lemma isogeny commutation}
    Let $x_1, x_2\in \mathfrak{M}^{\gr}(R)$, and suppose we are given an isogeny $f: \Ab{x_1}\to \Ab{x_2}$ between the abelian surfaces classified by those points, respecting the filtrations on the $p$-divisible groups $\Ab{x_i}[p^\infty]$. Then, for all $\zeta\in \mupi(R)$, we have that $f$ deforms to a moprhism $f': \Ab{\zeta^{m_2(f)} \cdot x_1}\to \Ab{\zeta^{m_0(f)} \cdot x_2}$, i.e.\ there exists a unique such $f'$ such that $f\equiv f' \mod (\zeta-1)$. Moreover, $f'$ is an isogeny that respects the filtrations on the $p$-divisible groups, and $m_i(f')=m_i(f)$ for all $i\in \lbrace 0,1,2\rbrace$.
    \begin{proof}
        By Serre-Tate lifting theory it is enough to show that the deformation $f'$ we are looking for can be built between the $p$-divisble groups of the two abelian schemes, i.e.\ that there exists $f': A[p^\infty]+E_{\zeta^{m_2(f)}}\to A[p^\infty]+E_{\zeta^{m_0(f)}}$ such that $f\equiv f' \mod (\zeta-1)$.

        Now, let $g_{m_0(f), m_2(f)}: E_{\zeta^{m_2(f)}}\to E_{\zeta^{m_0(f)}}$ be the morphism constructed in \Cref{lemma morphism Kummer}. Since $f$ and $g_{m_0(f), m_2(f)}$ have the same behavior once restricted to $\mupi$ and to $\QpZp$ (they both restrict to $m_0(f)$ and to $m_2(f)$, respectively), one can combine them to form a morphism $f': A[p^\infty]+E_{\zeta^{m_2(f)}}\to A[p^\infty]+E_{\zeta^{m_0(f)}}$ between the Baer sum of their respective domains and codomains, which will satisfy $m_i(f')=m_i(f)$ for all $i=0,1,2$. Modulo $\zeta-1$, $E_{\zeta^{m_0(f)}}$ and $E_{\zeta^{m_2(f)}}$ both coincide with the trivial extension $\mupi\oplus \QpZp$, $g_{m_0(f), m_2(f)}$ coincides with the morphism $m_0(f)\oplus m_2(f): \mupi\oplus \QpZp\to \mupi\oplus \QpZp$ by construction, and $f'$ consequently coincides with $f$.
    \end{proof}
\end{lemma}

We are now ready to state the following result about the interaction between the action of $\mupi$ and that of $\mathfrak{M}^\gr$.
\begin{proposition}
    \label{prop weight interaction mupi}
    Given $(a_0, a_1)\in \mathfrak{M}^{\gr}(R), \zeta\in \mupi(R)$ and a point $x\in \mathfrak{T}^\gr(R)$ away from the boundary, we have that $(a_0, a_1) \cdot \zeta \cdot (a_0^{-1}, a_1^{-1}) \cdot x = \zeta^{a_0^{-2}} \cdot x$.
    \begin{proof}
        As we noticed in \Cref{rmk m_i for id}, the identity morphism $f: \Ab{(a_0^{-1}, a_1^{-1})\cdot x}\to \Ab{x}$ is such that $m_2(f)=a_0$, $m_1(f)=a_1^{-1}$, and $m_0(f)=a_0^{-1}$. If we now apply \Cref{lemma isogeny commutation}, we get an isomorphism $f': \Ab{\zeta\cdot (a_0^{-1}, a_1^{-1})\cdot x}\to  \Ab{\zeta^{a_0^{-2}}\cdot x}$ such that $m_i(f')=m_i(f)$ for all $i\in \lbrace 0,1,2\rbrace$, which, in light of \Cref{rmk m_i for id}, is the same thing as an isomorphism $f'': \Ab{(a_0,a_1)\cdot \zeta \cdot (a_0^{-1}, a_1^{-1}) \cdot x}\to \Ab{\zeta^{a_0^{-2}}\cdot x}$ such that $m_i(f'')=1$ for all $i$.
    \end{proof}
\end{proposition}

\Cref{prop weight interaction mupi} can also be reformulated as follows.
\begin{corollary}
    \label{cor weight interaction contZp}
    Let $R$ be a $p$-adically complete $\Zp$-algebra, and $s$ a section of $\pi_*(\mathfrak{T}^\gr)$ over some open affine subscheme $U\inj \mathfrak{X}_{\mathrm{Kli}(p^\infty)}^{\ge 1}\cotimes_{\Zp}R$ that does not intersect the boundary $D$; moreover, let $f\in \Cont(\Zp,R)$ and $(a_0, a_1) \in \mathfrak{M}^{\gr}(R)$. Then, we have that 
    $(a_0, a_1) \cdot f \cdot (a_0^{-1}, a_1^{-1}) \cdot s = f^{[a_0]} \cdot s$, where $f^{[a_0]}(x):=f(a_0^{-2} x)$.
    \begin{proof}
        This is equivalent to \Cref{prop weight interaction mupi} if one re-interprets the action of $\mupi$ on $\mathfrak{T}^\gr$ as an action of $\Cont(\Zp, R)$ on $\pi_*(\mathfrak{T}^\gr)$, as explained in \S\ref{subsec action by continuous functions}.
    \end{proof}
\end{corollary}

If $f: \Zp \to R$ is a multiplicative function, given a character $k: \Zp^\times\to R$ one can form a new character $k+2f: \Zp^\times \to R$, $m\mapsto f(m)^{2} k(m)$. With this notatation, we can state the following result, which is a direct consequence of \Cref{cor weight interaction contZp}.
\begin{theorem}
    \label{thm main result}
    For any pair $(k_1,k_2)$ with $k: \Zp^\times\to R$ a character (for $R$ some $p$-adically complete $\Zp$-algebra) and $r\in \ZZ$, each multiplicative continuous function $f: \Zp\to R$ induces an $R$-linear operator $\vartheta^f: \mathfrak{F}^{(k_1,k_2)}\to \mathfrak{F}^{(k_1+2f,k_2)}$ over $\mathfrak{Y}_{\mathrm{Kli}(p^\infty)}^{\ge 1}$. In particular, for $n\in \NN$, the function $\Zp \to R$, $x\mapsto x^n$ induces a $n$-th order $R$-linear differential operator $\vartheta^n: \mathfrak{F}^{(k_1,k_2)}\to \mathfrak{F}^{(k_1+2n,k_2)}$.
\end{theorem}

\subsection{The effect on Hecke eigenvalues}
Let us now move to the interaction between the action of $\mupi$ and Hecke operators. Let us fix a matrix $w=\diag(d_1, d_2, d_3, d_4)\in \GSp_4(\mathbb{Z})$, and assume no $d_i$ is divisible by $p$. We necessarily have that $d_3 d_2 = d_4 d_1=:d(w)$. To such a coset, we can attach a correspondence $\mathfrak{C}\inj \mathfrak{U}\times_{\Zp} \mathfrak{U}$ that classifies the pairs of points $(x_1, x_2)\in \mathfrak{U}:= \mathfrak{T}^{\gr}\setminus D$ such that there exists a prime-to-$p$ isogeny $f: \Ab{x_1}\to \Ab{x_2}$ of type $w$, meaning that
\begin{enumerate}[(a)]
    \item the kernel $\ker(f)$ is a totally isotropic subgroup of $\Ab{x_1}[d(w)]$, and the isogeny $f$ respects the principal polarisations $\lambda_i$ of the abelian surfaces $A(x_i)$, meaning that the following diagram commutes
    \begin{equation*}
        \begin{tikzcd}
            A(x_1) \ar["f"]{r}\ar["d(w) \lambda_1"]{d}& A(x_2)\ar["\lambda_2",swap,"\cong"']{d}\\
            A(x_1)^\vee & A(x_2)^\vee\ar["f^\vee"]{l}
        \end{tikzcd}
    \end{equation*}
    \item for some appropriate choices of bases of the prime-to-$p$-Tate modules of $\Ab{x_1}$ and $\Ab{x_2}$, $f$ is represented by the matrix $w$;
    \item the isogeny $f$ respects the filtration on the $p$-disibible groups of $A(x_1)$ and $A(x_2)$, and its behaviour $m_i(f)$ on the trivialized graded piece $\gr_i(A[p^\infty])$ is fixed as follows: $m_0(f)=1, m_1(f)=1$ and $m_2(f)=d(w)$.
\end{enumerate}
Let $\pi_1, \pi_2: \mathfrak{C}\to \mathfrak{U}$ be the projections. It is clear that $\mathfrak{C}\inj \mathfrak{U}\times_{\Zp}\mathfrak{U}$ is invariant under the diagonal action of $\mathfrak{M}^\gr$, and that the maps $\pi_1, \pi_2: \mathfrak{C}\to \mathfrak{U}$ are both $\mathfrak{M}^\gr$-equivariant. The correspondence induces a map in cohomology $T_w: R\Gamma(\mathfrak{U};(k_1, k_2))\to R\Gamma(\mathfrak{U};(k_1, k_2))$, which is obtained by composing the pullback map $\pi_2^*: R\Gamma(\mathfrak{U};(k_1, k_2))\to R\Gamma(\mathfrak{C};(k_1, k_2))$ and the trace map $\mathrm{tr}(\pi_1): R\Gamma(\mathfrak{C};(k_1, k_2))\to R\Gamma(\mathfrak{U};(k_1, k_2))$. This is the desired Hecke operator for the forms of weight $(k_1, k_2)$.

Our aim now is studying the behaviour of the correspondence with respect to the action of $\mupi$. Let us consider the action by $\mupi$ on $\mathfrak{U}\times_{\Zp}\mathfrak{U}$ defined as follows: $\zeta\cdot (x_1, x_2):= (\zeta^{m_2(w)}\cdot x_1, \zeta^{m_0(w)}\cdot x_2)$. By \Cref{lemma isogeny commutation}, this action restricts to an action by $\mupi$ on  $\mathfrak{C}$. It is clear that, this time, $\pi_1$ and $\pi_2$ are not equivariant, but interact with the action in a non-trivial way:
\begin{equation*}
\begin{tikzcd}
    \mathfrak{C}\ar["\zeta"]{r} \ar["\pi_2"]{d} & \mathfrak{C}  \ar["\pi_2"]{d}\\
    \mathfrak{U}\ar["\zeta^{m_0(w)}"]{r} & \mathfrak{U}
\end{tikzcd}
\qquad
\begin{tikzcd}
    \mathfrak{C}\ar["\zeta"]{r} \ar["\pi_1"]{d} & \mathfrak{C}  \ar["\pi_1"]{d}\\
    \mathfrak{U}\ar["\zeta^{m_2(w)}"]{r} & \mathfrak{U}
\end{tikzcd}
\end{equation*}
as a consequence, we get the following result at the level of cohomology:
\begin{equation*}
\begin{tikzcd}
    R\Gamma(\mathfrak{U})\ar["\zeta^{{m_0(w)}}"]{d} \ar["\pi_2^*"]{r} & 
    R\Gamma(\mathfrak{C}) \ar["\zeta"]{d} \ar["\mathrm{tr}(\pi_1)"]{r} &
    R\Gamma(\mathfrak{U}) \ar["\zeta^{m_2(w)}"]{d} 
    \\
    R\Gamma(\mathfrak{U}) \ar["\pi_2^*"]{r} & 
    R\Gamma(\mathfrak{C}) \ar["\mathrm{tr}(\pi_1)"]{r} &
    R\Gamma(\mathfrak{U}) 
\end{tikzcd}
\end{equation*}

In other words, we have that $T_w \circ \zeta = \zeta^{m_2(w)/m_0(w)} \circ T_w$. Meanwhile, $m_2(w)/m_0(w) = d(w)/1 = d(w)$. If we rephrase this commutation relation in terms of actions by $\Cont(\Zp,R)$, we have proved the following proposition.
\begin{proposition}   \label{prop commutation with Hecke prime to p} Let $f: \Zp\to R$ be a multiplicative continuous function, and let $\vartheta^f: \mathfrak{F}^{(k_1,k_2)}\to \mathfrak{F}^{(k_1+2f,k_2)}$ be the corresponding differential operator. Then, $T_w \vartheta^f = f(d(w)) \vartheta^f T_w$.
\end{proposition}

This can in particular be applied to the operators $T_{\ell, 0}$, $T_{\ell, 1}$ and $T_{\ell,2}$ introduced in \S\ref{subsec hecke alg} to obtain commutation relations that are analogous to those presented in the introduction for the elliptic $\vartheta$ operator.
\begin{corollary}
    \label{cor commutation with Hecke prime to p}
    The following commutation relations hold, for all primes $\ell \nmid Np$.
    \begin{equation*}
    \begin{aligned}
    &T_{\ell,0} \vartheta = \ell^2 \vartheta T_{\ell,0}\\
    &T_{\ell,1} \vartheta = \ell^2 \vartheta T_{\ell,1}\\
    &T_{\ell,2} \vartheta = \ell \vartheta T_{\ell,2}
    \end{aligned}
    \end{equation*}
\end{corollary}

\section{Application to Galois representations}
\label{sec applications}
\subsection{Setup} For each weight $(k_1,k_2)\in \mathbb{Z}^2$, $k_1\ge k_2$, one can compute the coherent cohomology groups $H^i(X, \omega^{(k_1, k_2)})\otimes\Qp$ and $H^i(X, \omega^{(k_1, k_2)}(-D))\otimes \Qp$ of the Siegel threefold, that are finite-dimensional $\Qp$-vector spaces, endowed with an action of the Hecke algebra $\mathcal{H}$. Suppose now $f$ is a cohomology class belonging to one of these cohomology groups, and that it is a simultaneous Hecke eigenform for all Hecke operators; its eigenvalues can be packed together into a homomorphism $\Theta_f: \mathcal{H}\onto \Qpbar$.

Let now $\rho: G_\QQ\to \GSp_4(E)$ be a $p$-adic Galois representation (being $E$ some finite extension of $\Qp$) that is unramified away from the primes dividing $Np$, and let $\sigma_\ell\in G_\QQ$ be a Frobenious element relative to some prime $\ell\nmid Np$, i.e.\ $\sigma_\ell$ belongs to the decomposition group $D_{\mathfrak{L}}$ for some place $\mathfrak{L}$ of $\Qbar$ above $\ell$, and induces the Frobenious automorphism $x\mapsto x^\ell$ on the redidue field $k(\mathfrak{L})=\overline{\mathbb{F}_\ell}$. By unramifiedness, the conjugacy class of $\rho(\sigma_\ell)$ only depends on $\ell$; in particular, for each prime $\ell$ one can compute the characteristic polynomial $\chi_\ell(\rho)(x):=\det(1-x\rho(\sigma_\ell))$, which is a degree-4 polynomial with coefficients in $\Qpbar$ and constant term equal to 1. These polynomial characterize the representation uniquely, as a consequence of Cebotarev density theorem.
\begin{definition}
    We say that a semisimple continuous Galois representation $\rho: G_\QQ\to  \GSp_4(\Qpbar)$ unramified away from the primes dividing $Np$ \emph{is attached} to a Hecke eigensystem $\mathcal{H}^{Np}\onto \Qpbar$ if $\chi_{\ell}(\rho)=\Theta(Q_\ell)$ for all $\ell \nmid Np$, being $Q_\ell$ the Hecke polynomial as defined in \Cref{subsec hecke alg}.
\end{definition}

For Hecke eigensystems contributing to the coherent cohomology of the Siegel threefold $X$, we have the following existence result, which is a rephrasing of \cite[Theorem 5.3.1]{pilloniGSP4}.
\begin{theorem}
    \label{thm ext Galois}
    Let $\Theta: \mathcal{H}\to E\subseteq\Qpbar$ be a Hecke eigensystem such that $\Theta = \Theta_f$ for some eigenform $f\in H^i(X\otimes\Qpbar,\omega^{(k_1,k_2)})$ or $f\in H^i(X\otimes\Qpbar,\omega^{(k_1,k_2)}(-D))$. Then, there exists a semisimple Galois representation $\rho_f: G_\QQ\to  \GSp_4(E)$ attached to $\Theta$, unramified away from $Np$.
\end{theorem}

\subsection{Twisting by cyclotomic characters}
Let $\omega: G_\QQ\to \GL_1(\Qp)$ denote the cyclotomic character, characterized by the fact that it is unramified from $p$, and, for each prime $\ell\neq p$, $\omega(\ell)=\ell$. As a consequence of these properties, given a Galois representation $\rho: G_\QQ\to \GSp_4(\Qpbar)$, its cyclotomic twist $\rho\otimes \omega$ can be described as follows.
\begin{proposition}
    \label{properties cyclotomic twists}
    Let $\rho$ and $\rho\otimes \omega$ be as above. Then, given a prime $\ell\neq p$, $\rho\otimes \omega$ is unramified at $\ell$ if and only if $\rho$ is; moreover, in this case, $\chi_\ell(\rho\otimes \omega)(x)=\chi_\ell(\rho)(\ell x)$.
\end{proposition}

Applying the theta operator $\vartheta$ we have defined to a modular form is related to taking the cyclotomic twist of its attached Galois representation in the following way.
\begin{theorem}
\label{thm Galois}
    Let $f\in H^i(X,\omega^{(k_1,k_2)})\otimes \Qpbar$ be an eigenform for the Hecke algebra $\mathcal{H}$, being $i=0$ or $i=1$; let $\rho_f$ be its attached Galois representation. Let $g:=\vartheta(f)\in H^i(\mathfrak{Y}_{\mathrm{Kli}(p^\infty)}^{\ge 1},\mathfrak{F}^{(k_1+2,k_2)})\otimes \Qpbar$. Assume $g\neq 0$. Then, $g$ is still an eigenform for the prime-to-$p$ Hecke algebra $\mathcal{H}^{p}$; moreover, there exists a Galois representation $\rho_g$ attached to its Hecke eigensystem, which coincides with the cyclotomic twist of $f$ (in other words, $\rho_g\cong \rho_f\otimes \omega$).
    \begin{proof}
        It follows from the commutation properties between Hecke operators and $\vartheta$ (\Cref{cor commutation with Hecke prime to p}) that, if $a_{\ell,i}$ (being $i=0,1,2$) is the Hecke eigenvalue of $f$ with respect to $T_{\ell,i}$, then $g$ is an eigenform for $T_{\ell,i}$ with respect to the eigenvalue $\ell^{M(i)} a_{\ell_i}$, being $M(0)=M(1)=2$, and $M(2)=1$.  Let us denote by $\Theta_f$ and $\Theta_g$ the Hecke eigensystems of $f$ and $g$ respectively: recalling the homogeneity properties of the Hecke polynomial $Q_\ell$ described in \Cref{rmk Ql homogeneous}, we consequently have that $\Theta_g(Q_\ell)(x) = \Theta_f(Q_\ell)(\ell x)$.

        The representation $\rho_f$ is unramified at all $\ell\nmid Np$, and the characteristic polynomials of the Frobenii are $\chi_\ell(\rho_f)(x)=\Theta_f(Q_\ell)(x)$ (see \Cref{thm ext Galois}). If we recall the description of the cyclotomic twist given in \Cref{properties cyclotomic twists}, we clearly have that $\rho_f\otimes \omega$ is also unramified at all $\ell\nmid Np$, and that $\chi_\ell(\rho_f\otimes \omega)(x)=\Theta_f(Q_\ell)(\ell x)$.
        
        Putting everything together, we conclude that $\chi_\ell(\rho_f\otimes \omega)(x)=\Theta_g(Q_\ell)(x)$, which is to say that $\rho_f\otimes \omega$ is the Galois representation attached to the Hecke eigensystem of $g$.
    \end{proof}
\end{theorem}

\printbibliography

%If $c\in H^d(\mathfrak{X}, R)$ is a cohomology class (for $R$ a p-adically complete $\Zp$-algebra), in particular, we have that $(m_1, m_0) \cdot \zeta \cdot (m_1^{-1}, m_0^{-1}) \cdot x = \zeta^{m_0^2} \cdot x$ for all $(m_1, m_0)\in \mathfrak{M}^{\gr}(R)$ and $\zeta\in \mupi(R)$.

%\subsection{$p$-depletion}

%\input{general thing}
\end{document}